\documentclass[12pt,a4paper]{article}
\usepackage{amssymb}
\usepackage{amsmath}
\usepackage[dvips]{graphics}
\usepackage{graphicx}
\usepackage{latexsym}
\usepackage[normal]{subfigure}
\input{epsf.sty}
\include{psfig}
\usepackage{amsfonts}
\begin{document}
\def\K{\mathbb{K}}
\def\R{\mathbb{R}}
\def\C{\mathbb{C}}
\def\Z{\mathbb{Z}}
\def\Q{\mathbb{Q}}
\def\D{\mathbb{D}}
\def\N{\mathbb{N}}
\def\T{\mathbb{T}}
\def\P{\mathbb{P}}
\renewcommand{\theequation}{\thesection.\arabic{equation}}
\newtheorem{theorem}{Theorem}[section]
\newtheorem{cond}{C}
\newtheorem{lemma}{Lemma}[section]
\newtheorem{corollary}{Corollary}[section]
\newtheorem{prop}{Proposition}[section]
\newtheorem{definition}{Definition}[section]
\newtheorem{remark}{Remark}[section]

\title{Geometric generalizations in Kresin-Maz'ya
Sharp Real-Part Theorems \thanks{AMS classification number:30A10.}
\thanks{Keywords: Sharp constants, Real Part Theorems.}}
\author {Lev Aizenberg \\
 Department of Mathematics, Bar-Ilan University,\\
 52900 Ramat-Gan, Israel,
  email: aizenbrg@math.biu.ac.il \\
  and\\
 Alekos Vidras,\\
  Department of Mathematics and Statistics,\\
  Univ.of Cyprus, Nicosia 1678, Cyprus, \\email: msvidras@ucy.ac.cy}

\maketitle
 \begin{abstract}
 In the present article we give geometric generalizations
 of the estimates from Chapters 5,6,7 from \cite{krem:gnus},
 while extending their sharpness to new cases.
 \end{abstract}
 \section{Preliminaries }
\setcounter{equation}{0}
  G.Kresin and V.Maz'ya, in their recently published, remarkable
research monograph \cite{krem:gnus}, have collected in one place
generalizations and different modifications of theorems, which the
authors called {\it{real-part theorems}} honoring the known
theorem of Hadamard (1892). All of them are formulated with sharp
constants.\par
 In the present article, our starting point is the content of the
 three chapters of the above monograph, namely, Chap.5,{\it{ Estimates for the
derivatives of analytic functions,}} Chap.6, {\it{Bohr's type real
estimates}}, Chap.7, {\it{Estimates for the increment of
derivatives of analytic functions}}. Methods, used in
\cite{aiz:gnus}, allow a geometric generalization for some of the
results in those chapters
 and the sharpness of the results is extended to new cases. We
 remark, that there are two methods to prove the sharpness of the
 corresponding estimates. The first one, which is shorter, uses
 a number of facts from the monograph \cite{krem:gnus}, the
 second, somewhat longer, but independent of \cite{krem:gnus}. To
 illustrate the point we use the approach in \S2, and the second
 one in \S3.
\section{Estimates for the derivatives of holomorphic functions}
\setcounter{equation}{0} We begin by formulating a theorem,
inspired by Th.5.1, from \cite{krem:gnus}

\begin{theorem}
{\rm{(G.Kresin-V.Maz'ya)}}
 Let $f$ be holomorphic in
$\mathcal{D}_R=\{z\in\mathbb{C}:\; \vert z\vert <R\}$ and let
$z=re^{i\theta}$, $a=r_ae^{i\theta}$, $0\leq  r_a\leq r<R$. Then
the inequality
\begin{eqnarray}
\vert f^{(n)}(z)\vert \leq
{2n!R(R-r_a)\over(R-r)^{n+1}(R+r_a)}\mathcal{Q}_a(f)
\end{eqnarray}
holds for every $n\geq 1$ with the best possible constant , where
$\mathcal{Q}_a(f) $ is each of the following expressions\\
{\rm{i)} }$\sup_{\vert \zeta \vert <R}\Re f(\zeta)-\Re f(a)$. In
this case the claim of the {\rm Theorem 2.1} is a generalization
of the
Hadamard real-part theorem.\\
{\rm{ii)}} $\sup_{\vert \zeta \vert <R}\vert \Re f(\zeta)\vert
-\vert \Re f(a)\vert $. In this case the claim of the {\rm Theorem
2.1} is a
generalization of the Landau type inequality.\\
{\rm ii)} $\sup_{\vert \zeta \vert <R}\vert f(\zeta)\vert -\vert
f(a)\vert $. In this case the claim of the {\rm Theorem 2.1} is a
generalization of the Landau  inequality.\\
{\rm iv)} $\Re f(a)$, if $\Re f>0$ on $\mathcal{D}_R$. In this
case the claim of the {\rm Theorem 2.1} is a generalization of the
Caratheodory inequality.
 \end{theorem}
 We remark that the estimate $(2.1)$ has
a meaning in the cases {\rm(i)}, {\rm(ii)}, {\rm(iii)} only if the
corresponding $sup$ is finite. In the case $a=z$ the inequalities
are due to S.Ruscheweyh \cite{rush:gnus}. Other references on
different type of inequalities preceding  the inequality $(2.1)$
are to be found in \cite{krem:gnus}.\par
 Denote now by $\widetilde G$ the convex
hull of the domain $G\subset\mathbb{C}$. A point $p\in\partial G$
is called a {\em point of convexity } if $p\in \partial \widetilde
G$. A point of convexity $p$ is called {\em regular} if there
exists a disk $\mathcal {D}^{\prime}\subset G$ so that $p\in
\mathcal {D}^{\prime}$.
\begin{theorem}
Let $f$ be holomorphic function in $\mathcal{D}_R$,
$f(\mathcal{D}_R)\subset G$, where $G$ is a domain $\mathbb{C}$, so
that $\widetilde G\not=\mathbb{C}$. Let also $z=re^{i\theta}$,
$a=r_ae^{i\theta}$ be complex numbers so that $0\leq r_a\leq r<R$.
Then the inequality
\begin{eqnarray}
\vert f^{(n)}(z)\vert \leq {2n!R(R-r_a)\over
(R-r)^{n+1}(R+r_a)}dist(f(a),\partial \widetilde G)
\end{eqnarray}
holds for $n\geq 1$. If $\partial  G$ contain at least one regular
point of convexity, then the constant in $(2.2)$ is sharp.
\end{theorem}
\begin{remark}
We point out that the cases {\rm(i)} and {\rm(ii)} in
 {\rm Theorem 2.1} are the cases when
\begin{eqnarray*}
G=\{z\in\mathbb{C}:\; \Re z<\sup_{\vert \zeta\vert <R}\Re f(\zeta)\}
\end{eqnarray*}
 is half-plane or
 \begin{eqnarray*}
G=\{z\in\mathbb{C}:\; \vert z\vert <\sup_{\vert \zeta\vert <R}\vert
\Re f(\zeta)\vert\}
\end{eqnarray*}
 is a strip. The case {\rm(iii)} in {\rm Theorem 2.1} corresponds to the set
\begin{eqnarray*}
G=\{z\in\mathbb{C}:\; \vert z\vert <\sup_{\vert \zeta\vert <R}\vert
f(\zeta)\vert\}
\end{eqnarray*}
being a disc. The case {\rm(iv)} of this theorem corresponds to
the right half-plane case
\begin{eqnarray*}
G=\Pi=\{z\in\mathbb{C}:\; \Re f( z)>0\}.
\end{eqnarray*}
\end{remark}
{\bf{Proof:}} Let us consider the case ${\rm (iv)}$ in Theorem
2.1. If $\Re f(z)>0$, $z $ on $\mathcal{D}_R$, then for every
$n\geq 1$ one has
\begin{eqnarray}
\vert f^{(n)}(z)\vert \leq {2n!R(R-r_a)\over
(R-r)^{n+1}(R+r_a)}dist(f(a),\partial \Pi),
\end{eqnarray}
where $\Pi $ is the right half-plane. By translation and rotation,
this half plane can be transformed into any half-plane $\Pi _1$.
The same transformations can be applied to holomorphic in
$\mathcal{D}_R$ functions, that is $f(z)\longrightarrow
(f(z)+c)e^{i\phi }$. Under such transformations of both, the
half-plane and the functions, at the same time the conclusion
$(2.3)$ does not alter.\par
 Now, let $G$ be a domain in
$\mathbb{C}$, such that $\widetilde G\not=\mathbb{C}$. Let also
the distance from the right hand-side of the inequality $(2.2)$ is
realized at the point $p\in\partial\widetilde G$. Then there exist
a line of support to $\widetilde G$ at $p$, bounding the
half-plane $\Pi_1 $, such that $G\subset \Pi_1$. For this
particular half-plane $\Pi _1$ we apply the inequality $(2.3)$ and
using the fact
\begin{eqnarray*}
dist(f(a),\partial \widetilde G)=dist(f(a),p)=dist(f(a),\Pi_1),
\end{eqnarray*}
we obtain $(2.2)$. \par
 Assume now that $\partial G$ contains at
least one regular point of convexity $p\in \partial G\cap\partial
\widetilde G\cap\partial \mathcal{D}^{\prime}$, where
$\mathcal{D}^{\prime}$ is some disc $\mathcal{D}^{\prime}\subset
G$. Let $\mathcal{D}_{\beta}$ be a disk, in which the functions
from {\rm(iii)} of Theorem 2.1 take their values:
\begin{eqnarray*} \mathcal{D}_{\beta}=\{z\in\mathbb{C}:\vert
z\vert <\beta =\sup_{\vert \zeta \vert <R}\vert f(\zeta)\vert \}.
\end{eqnarray*}
For this disc $\mathcal{D}_{\beta}$ the constant in $(2.1)$ is
sharp. That is, there exists a family of functions, which we take
from \S5.7,\cite{krem:gnus},
\begin{eqnarray}
g_{\xi}(z)={\xi \over z-\xi}+{\vert \xi\vert^2\over \vert \xi
\vert ^2-R^2},
\end{eqnarray}
depending on complex parameter $\xi=\rho e^{i\theta}$, $\rho >R$,
for which this constant is attained. Put $z=x<0$, $\xi =\rho >0 $,
the sharpness of the constant in \cite{krem:gnus} was proved by
passing to the limit when $\rho\downarrow R$. Then the first
summand in $(2.4)$ is negative and remains bounded  while $\rho
 \downarrow R$. The second summand in $(2.4)$ tends to $+\infty $ while $\rho
\longrightarrow R$. Hence, when $\rho $ is sufficiently close to
$R$, the function $g_{\rho}(x)$ is positive. This implies that
$g_{\rho}(a)$ is also positive.\par Of crucial importance in our
reasoning is the fact that one can transform the disk
$\mathcal{D}_{\beta}$ into the disc $\mathcal{D}^{\prime}$ by
using homothety and translation. At the same time , we apply both
transformations to the family of functions
$g_{\xi}(z)\longrightarrow \widetilde g_{\xi}(z)=\alpha
g_{\xi}(z)+c$. Then the inequality $(2.1)$ in the case {\rm (iii)}
of the Theorem 2.1 does not change. Furthermore, if the for the
family of functions $\{g_{\xi}(z)\}_{\xi }$ the sharpness of the
constant was attained before the applications of the
transformations under the assumption $f(\mathcal{D}_R)\subset
\mathcal{D}_{\beta}$, then the family of functions obtained after
the transformation illustrates the sharpness of the constant under
the assumption $f(\mathcal{D}_R)\subset \mathcal{D}^{\prime}$.\par
 In the \S5.7 of \cite{krem:gnus}, for the proof of the sharpness of the given constant only the modula of $\vert
 g_{\xi}(z)\vert $ and  $\vert g_{\xi}(a)\vert $ were used.
 Therefore, instead of the family $\{g_{\rho}(z)\}_{\rho }$ one can
 use the family $\{e^{i\phi }g_{\rho}(z)\}_{\rho }$ in order to
 show the required sharpness of the constant, provided that $\phi $ is chosen in a such a manner  that the point
 $e^{i\phi }g_{\rho}(a)$ on the radius eminnating  from the center of the disc $ \mathcal{D}_{\beta}$ to
 the point $p^{\prime}\in\partial \mathcal{D}_{\beta}$, where $p^{\prime}$ is the pre-image of the point $p$ under the
 above mentioned homothety and translation  of the disc  $ \mathcal{D}_{\beta}$ into the disc  $\mathcal{D}^{\prime}$.
Then, after the transformation, the point $\widetilde g_{\rho}(a)$
will lie on the radius eminnating from the center of the disc
$\mathcal{D}^{\prime}$ to the point $p$, and therefore
\begin{eqnarray*}
dist(\widetilde g_{\rho}(a),\partial \widetilde G)=dist(\widetilde
g_{\rho}(a),p).
\end{eqnarray*}
This completes the proof of the theorem. $\diamondsuit $\par
 For $a=0$ one has the following
 \begin{corollary}
 Let $f$ be holomorphic in the disc $ \mathcal{D}_{R}$ and $f(\mathcal{D}_{R})\subset
 G$, where $G$ is a domain in $\mathbb{C}$ whose convex hull
 $\widetilde G $ is not equal to $\mathbb{C}$.Then the inequality
 \begin{eqnarray}
 \vert f^{(n)}(z)\vert\leq {2n!R\over (R-r)^{n+1}}dist (f(0),
 \partial \widetilde G)
 \end{eqnarray}
 holds for every $n\geq 1$.
 If $\partial \widetilde G $ contains at least one regular point
 of convexity , then the constant in $(2.5)$ is sharp.
 \end{corollary}
 \section{Bohr's type real part estimates}
\setcounter{equation}{0} The results, contained in the Theorems
6.1-6.4 in \cite{krem:gnus}, are collected in the following
\begin{theorem}
({\rm G.Kresin-V.Maz'ya}) Let the function
\begin{eqnarray}
f(z)=\sum\limits_{n=0}^{\infty}c_nz^n
\end{eqnarray}
be holomorphic in the disc $ \mathcal{D}_{R}$ and $q>0$, $m\geq
1$, $\vert z\vert =r <R $. Then the inequality
\begin{eqnarray}
\left(\sum\limits_{n=m}^{\infty}\vert c_nz^n\vert
^q\right)^{{1\over q}}\leq {2r^m\over R^{m-1}(R^q-r^q)^{{1\over
q}}}\mathcal{R}(f)
\end{eqnarray}
holds with the best possible constant in the cases when
$\mathcal{R}(f)$ is each one of the following expressions:\\
{\rm i)} $\sup_{\vert \zeta \vert <R}(\Re (f(\zeta)-\Re f(0))$.\\
{\rm ii)} $\sup_{\vert \zeta \vert <R}(\vert \Re f(\zeta)\vert
-\vert
\Re f(0)\vert )$. \\
{\rm iii)} $\sup_{\vert \zeta \vert <R}(\vert f(\zeta)\vert -\vert
f(0)\vert )$. \\
{\rm iv)} $\Re f(0)$, if $\Re f>0$ on $\mathcal{D}_R$.
\end{theorem}
We remark here that the case ${\rm(iii)}$ of the above theorem
gives for $m=q=1$ the classical theorem of Bohr, \cite{bohr:gnus}
(for related references see \cite{aiz:gnus}) for $r={R\over 3}$.
Similarly to the previous section, we state the geometric
generalization of this theorem.
\begin{theorem}
Let the function $(3.1)$ be holomorphic in the disc $\mathcal{D}_R$,
$q>0$, $m\geq 1$, $\vert z\vert =r<R$, and $f(\mathcal{D}_R)\subset
G$, where $G$ be a domain in $\mathbb{C}$ and $\widetilde
G\not=\mathbb{C}$. Then the inequality
\begin{eqnarray}
\left(\sum\limits_{n=m}^{\infty}\vert c_nz^n\vert
^q\right)^{{1\over q}}\leq {2r^m\over R^{m-1}(R^q-r^q)^{{1\over
q}}}dist(c_0,\partial \widetilde G),
\end{eqnarray}
holds. If the boundary $\partial G$ contains at least one regular
point of convexity, then the constant in $(3.3)$ is sharp.
\end{theorem}
{\bf{Proof:}} The estimate $(3.3)$ is proven in exactly the same
way as in the previous section. The sharpness of the constant in
$(3.3)$ can also be proven in the same manner, but we prefer to
give an independent proof.\par The main point of our approach is
the simple observation that convex hull of the domain between to
discs, when the smaller one is contained in the larger one and
their boundaries have exactly one common point, is the larger
disc.\par
 To be more specific, for $a>0$ we consider
\begin{eqnarray*}
D_1&=& \{z\in {\bf C}: \vert z-ai\vert < a\}\\
D_2&=&\{z\in {\bf C}: \vert z-2ai\vert < 2a\}\\
\end{eqnarray*}
be two discs. It is clear that $\partial D_1\cap \partial
D_2=\{0\}$. Define the domain
\begin{eqnarray*}
G=D_1^c\cap D_2,
\end{eqnarray*}
where, as usual, $D_1^c $ denotes the complement of the disc $D_1$
in $\mathbb{C}$.\par
 It is obvious that the convex hull
$\widetilde G$ of $G$ is equal to $D_2$. The conformal map
$f(\zeta)={1\over \zeta }$ maps $G$ onto a strip
\begin{eqnarray*}
T=\{z\in {\bf C}: -{i\over 2a}<\Im z<-{i\over 4a}\}
\end{eqnarray*}

The $\partial T$ consists of two parallel lines, on which $z$ moves
in the opposite directions. The point $z=\infty $ is a double point.
The width of the strip is equal to ${1\over 4a}$. The map
$w=f_1(z)=z+{i\over 2a}$ shifts the strip up. The new strip $
T+{i\over 2a}$ has the real axis as its lower bound, while the width
remains the same. Then  the map $\omega=f_2(w)=e^{4a\pi w}$
transforms con-formally the strip onto the upper half plane
\begin{eqnarray*}
H=\{\omega \in {\bf C}: \Im \omega >0\}
\end{eqnarray*}
Finally, we transfer $H$ by translation $\psi=f_3(\omega)=\omega
+p$, where $p\in {\bf C}$, $\Im p <0 $. All the above maps are
invertible. Thus we have a map
\begin{eqnarray*}
F&:& H+p\longrightarrow G\\
F(\psi)&=&{1\over {1\over 4a\pi }\ln (\psi-p)-{i\over 2a}},\\
\end{eqnarray*}
which is the inverse of the composition of $f_i$. Its expansion in
the disc $D(0,\vert p\vert)$ is given by
\begin{eqnarray*}
c_n&=&{a_n-c_0b_n-c_1b_{n-1}-\dots-c_{n-1}b_1\over
b_0},\;\rm{where}\\
 a_n&=&0,\; \forall n\geq 1,\;a_0=1,\\
b_n&=&{({1\over 4a\pi}\ln(\psi
-p)-{i\over 2a})^{(n)}\over n!}\vert _{\psi =0},\;n\geq 1,\\
\end{eqnarray*}
 and $c_0=F(0)$. Or equivalently,
\begin{eqnarray*}
c_n&=&{(-1)^nb_1^n\over b_0^{n+1}}, \;\rm{where}\\
b_1&=&{1\over 4a\pi}{1\over (-p)},\;b_0={1\over
4a\pi}\ln(-p)-{i\over 2a}\\
\end{eqnarray*}

Thus
\begin{eqnarray*}
c_n={4a\pi\over p^n (\ln (-p)-{i\over 2a})},\;n\geq 1
\end{eqnarray*}
If $p=-i$ then $c_0={1\over {1\over 4a\pi}\ln (i) -{i\over
2a}}={8a\over 3}i$
 and hence
\begin{eqnarray*}
dist(c_0,\tilde G)={4a\over 3}
\end{eqnarray*}

Assume now, that the best constant in $(3.3)$ is denoted by $C(r)$.
For $\vert z\vert =r $, $0<r<p$, one has
\begin{eqnarray*}
\sum\limits_{n=1}^{\infty}\vert c_n z^n\vert
^q&=&\sum\limits_{n=1}^{\infty}({4\pi a\over \vert \ln(-p)-{i\over
2a}\vert})^q{r^{qn}\over p^{qn}}
\end{eqnarray*}

The last power series is geometric one with ratio $({r\over p})^q$.
Therefore, for $m\geq 1$, one has
\begin{eqnarray*}
({4a\pi \over \vert \ln(-p)-{i\over
2a}\vert})^q\sum\limits_{n=m}^{\infty}({r\over p})^{nq}&=&({4a\pi
\over \vert \ln(-p)-{i\over 2a}\vert})^q{r^{mq}p^q\over
p^{mq}(p^q-r^q)}\\
&=&2{4a\over 3}{r^{mq}p^q\over p^{mq}(p^q-r^q)}
\end{eqnarray*}

Thus, the estimate from above is
\begin{eqnarray*}
{8a\over 3} {r^{m}\over p^{m-1}}{1\over ((p)^q-r^q)^{{1\over
q}}}\leq C(r) dist(c_0,\tilde G)=C(r){4a\over 3}
\end{eqnarray*}
or
\begin{eqnarray*}
{2r^{m}\over  R^{m-1}((R^q-r^q))^{{1\over q}}}\leq C(r)
\end{eqnarray*}

Thus the exactness is proven for $m=1,2,\dots $. \par
This
particular example shows that in the case of bounded, simply
connected domain $G$, whose boundary contains a point of regular
convexity, the constant in the Theorem 3.2 is sharp. Actually if
$\zeta _0$ is a point of regular convexity, then it means that
$\zeta_0\in
\partial \widetilde {\mathcal G}$ and there is a disc $U$ of radius
$\rho$, contained in $ G$ so that $\partial U\cap\partial
{\mathcal G}=\{\zeta _0\}$. We inscribe in the disc $U$, the disc
$U_1$ whose center lies in the diameter of $U$, whose end-points
are $\zeta _0, \zeta_0^{\prime}$. The radius of the disc $U_1$ is
${1\over 2}\rho$ and its center $k_1$ satisfies $\vert
k_1-\zeta_0^{\prime}\vert ={\rho\over 2}$. Then for the domain
${\mathcal U}=U\cap U_1^c $ we repeat the construction of the
above example. The key fact here is that $\widetilde {\mathcal
U}=U $ and that the distance $d(c_0, \widetilde {\mathcal U})$ is
realized at the point $\zeta _0\in \widetilde G$. So the sharpness
of the constant is proven.$\diamondsuit $\par
 Furthermore, the Theorem 6.5, \cite{krem:gnus}, states
\begin{theorem}
({\rm G.Kresin-V.Maz'ya})
 Let $f(z)$ be a function holomorphic in
the disc $\mathcal {D}_R$ and assume that in the neighborhood of
the point $a \in \mathcal {D}_R$ the expansion
\begin{eqnarray}
f(z)=\sum\limits_{k=0}^{\infty}c_k(a)(z-a)^k
\end{eqnarray}
is valid. Then for every $z\in \mathcal {D}_R$, $\vert z-a\vert =r
<d_a =dist(a, \partial \mathcal {D}_R)$ the following inequality
\begin{eqnarray*}
\sum\limits_{k=1}^{\infty}\vert c_k(a)(z-a)^k\vert \leq {2Rr\over
(2R-d_a)(d_a-r)}\mathcal{Q}_a(f)
\end{eqnarray*}
holds with the best possible constant and where $\mathcal{Q}_a(f)$
is each of the following expressions {\rm(i)}, {\rm(ii)},
{\rm(iii)}, {\rm(iv)} from the {\rm Theorem 2.1}.
\end{theorem}
Similarly one can prove
\begin{theorem}
Let $f(z)$ be a function holomorphic in the disc $\mathcal {D}_R$
and assume that in the neighborhood of the point $a \in \mathcal
{D}_R$ the expansion $(3.4)$ is valid. Assume also that
$f(\mathcal {D}_R)\subset G$, where $G$ is a domain in
$\mathbb{C}$ such that $\widetilde G\not=\mathbb{C}$. Then for
every $z\in \mathcal {D}_R$, $\vert z-a\vert =r <d_a =dist(a,
\partial \mathcal {D})$ the following inequality
\begin{eqnarray}
\sum\limits_{k=1}^{\infty}\vert c_k(a)(z-a)^k\vert \leq {2Rr\over
(2R-d_a)(d_a-r)}dist(f(a), \partial \widetilde G)
\end{eqnarray}
holds. If $\partial G $ contains at least one regular point of
convexity, then the constant in $(3.5)$ is sharp.
\end{theorem}
\section{Estimates for the increment of derivatives of holomorphic functions}
\setcounter{equation}{0} In this section we will be using the
notation $\Delta g(z)=g(z)-g(0)$ to describe the increment of a
function $g$ at $z=0$. We will formulate the results of the
Corollaries 7.2-7.5 from the book \cite{krem:gnus}in the following
manner.
\begin{theorem}
({\rm G.Kresin-V.Maz'ya}) Let $f(z)$ be a function holomorphic in
the disc $\mathcal {D}_R$. Then, for any fixed $z$, $\vert z\vert
=r<R$, the inequality
\begin{eqnarray*}
\vert f^{(n)}(z)-f^{(n)}(0)\vert \leq
{2n!(R^{n+1}-(R-r)^{n+1})\over (R-r)^{n+1}R^n}\mathcal{R}(f)
\end{eqnarray*}
holds with the best constant for every $n\geq 0$ and where
$\mathcal{R}(f)$ is each of the expression {\rm(i)}-{\rm (iv)}
from the {\rm Theorem 3.1}.
\end{theorem}
Analogously to the previous sections one can prove the following
\begin{theorem}
Let $f(z)$ be a function holomorphic in the disc $\mathcal {D}_R$.
  Assume also that $f(\mathcal {D}_R)\subset G$, where $G$ is a
domain in $\mathbb{C}$ such that $\widetilde G\not=\mathbb{C}$.
Then for every fixed $z\in \mathcal {D}_R$, $\vert z\vert =r $ the
following inequality
\begin{eqnarray}
\vert f^{(n)}(z)-f^{(n)}(0)\vert \leq
{2n!(R^{n+1}-(R-r)^{n+1})\over (R-r)^{n+1}R^n}dist
(f(0),\partial\widetilde G)
\end{eqnarray}
holds for every $n\geq 0$. If $\partial G $ contains at least one
regular point of convexity, then the constant in $(4.1)$ is the
best one.
\end{theorem}
We conclude the article with the following
\begin{remark}
 1) It seems that it is possible to formulate a
geometric
variant of the results from \cite{kres:gnus}. \\
2) Do the constants in the above cited geometric generalizations
of results of Kresin-Maz'ya  remain sharp, if one assumes
 that $\partial G$ does not contain any  regular point of convexity?
\end{remark}

\end{document}